\documentclass[a4paper,11pt]{article}
\pdfoutput=1
\usepackage[a4paper,margin=30mm]{geometry}
\usepackage[sec]{style_unfolding}

\title{On closed symmetric monoidal enriched $\infty$-categories}
\author{Rune Haugseng}

\newcommand{\SMENRlc}{\SMENR^{\lax,\cl}_{\CA}}
\newcommand{\SMENRcl}{\SMENR^{\cl}}
\newcommand{\ARoplcl}{\AR^{\oplax}_{\mathrm{s}=\cl}(\SMCATl)}
\newcommand{\FFr}{\mathfrak{F}}
\newcommand{\tMor}{\tcatname{M}{or}}

\newcommand{\SM}{SM}
\newcommand{\LSM}{LSM}
\newcommand{\CSM}{CSM}
\date{\today}
\begin{document}

\maketitle

\begin{abstract}
  We show that a closed symmetric monoidal $\oV$-enriched
  $\infty$-category can be described by a lax symmetric monoidal
  functor to $\oV$ from its underlying \icat{}.
\end{abstract}

\tableofcontents

\section{Introduction}
Suppose $\eC$ is a symmetric monoidal (SM) category that is enriched in
$\oV$. If $\oC$ denotes the underlying ordinary category of $\eC$,
taking enriched Hom out of the monoidal unit gives a lax symmetric
monoidal (LSM) functor
\[ \phi = \eC(\bbone,\blank) \colon \oC \to \oV.\]
If the symmetric monoidal structure on $\oC$ is closed, \ie{} it has
an internal Hom, and this is compatible with the enrichment in the
sense that the canonical maps
\[ \eC(x, [y,z]) \to \eC(x \otimes y, [y,z] \otimes y) \to \eC(x
  \otimes y, z)\]
are isomorphisms, then we can reconstruct $\eC$ from the \LSM{} functor $\phi$: Since $\oC$ is closed symmetric
monoidal (CSM), it has an enrichment $\bar{\oC}$ in itself given by
the internal Homs, and we can push
forward the enrichment along $\phi$  to get a \SM{}
$\oV$-category with enriched Homs given by
\[ \phi([x,y]) = \eC(\bbone,[x,y]) \cong \eC(x,y),\]
so that $\phi_{*}\bar{\oC} \simeq \eC$.

Our goal in this paper is to prove that this description of \CSM{} enrichments remains true for \icats{}. More
precisely, we establish an equivalence between the following
\itcats{}:
\begin{itemize}
\item $\SMCATlcl_{\upslash \oV}$ is the \itcat{} of \LSM{}
  functors $\oC \to \oV$ where $\oC$ is closed, with morphisms oplax
  triangles of \LSM{} functors
    \[
    \begin{tikzcd}[column sep=small]
     \oC \ar[rr] \ar[dr, ""{name=A}] & & \oD \ar[dl] \\
     & \oV.
     \arrow[from=A,to=1-3,Rightarrow,shorten <=4pt, shorten >=10pt]
    \end{tikzcd}
  \]
\item $\PSMENRlcl(\oV)$ is the \itcat{} of \emph{pinned} \CSM{} $\oV$-\icats{} and \LSM{}
  functors. By this we mean that the objects are pairs $(\eC, \oC' \xto{f} U\eC)$
  where $\eC$ is a \CSM{} $\oV$-\icat{} and $f$ is
  an essentially surjective strong \SM{} functor that
  preserves internal Homs.
\end{itemize}
\begin{theoremint}\label{thm:main V}
  For a  \SM{} \icat{} $\oV$, there is an equivalence of \itcats{}
  \[ \SMCATlcl_{\upslash \oV} \simeq \PSMENRlcl(\oV). \]
  This takes a \LSM{} functor $\phi \colon \oC \to
  \oV$ to $\phi_{*}\bar{\oC}$ together with a canonical morphism $\oC
  \to U\phi_{*}\bar{\oC}$ to its underlying \icat{}. In the other
  direction a pair $(\eC, f \colon \oC' \to U\eC)$ corresponds to the
  composite
  \[ \oC' \xto{f} U\eC \xto{\eC(\bbone,\blank)} \oV.\]
  This restricts to an equivalence between the full sub-\itcat{} 
  $\SMENRlcl(\oV)$ of \CSM{} $\oV$-\icats{} and
  that of \LSM{} functors $\phi$ for which the
  canonical maps
  \[ \oC(\bbone,x) \to \oV(\bbone, \phi x)\]
  are equivalences.
\end{theoremint}
  
\begin{warning*}
  Our proof of \cref{thm:main V} makes use of a new model for
  $\oV$-enriched \iopds{}, and we assume without proof that this is
  2-functorial in $\oV$ by changing enrichment along \LSM{} functors
  in the expected way. See \cref{assumption:eopdfib} for a more
  precise statement and \cref{rmk:assumption sketch} for a sketch of a
  proof, which we expect to write up in more detail elsewhere.
\end{warning*}  

\begin{warning*}
  In \cref{thm:main V}, the terms \SM{} $\oV$-\icats{} and \LSM{}
  functors refer to a certain full subcategory of our model for
  $\oV$-\iopds{}. It has currently not been shown for \emph{any}
  definition of enriched \iopds{} that the usual definition of \SM{}
  $\oV$-\icats{} with \LSM{} functors (defined as commutative algebras
  in $\oV$-\icats{} and certain inert-lax transformations among them)
  gives a full subcategory thereof. However, in the special case of
  enrichment in $(\infty,n-1)$-categories that is relevant to
  \cite{unfolding}, such a description of symmetric monoidal
  $(\infty,n)$-categories in terms of the dendroidal model of enriched
  \iopds{} from \cite{enropd} will be proved in forthcoming work of
  Tallak Manum; this still leaves the issue of comparing this model to
  the one used here, however.
\end{warning*}

By applying the equivalence of \cref{thm:main V} on
mapping \icats{}
we obtain the following universal property\footnote{Note that this immediately implies the universal property of
\cite[Proposition 1.3.3]{GRV}, which is stated without a complete
proof there.} of pushforwards
$\phi_{*}\bar{\oC}$ :
\begin{corint}
  Let $\oC$ be a \CSM{} \icat{}, $\phi \colon \oC \to \oV$ a \LSM{}
  functor, and $\eD$ a \CSM{} $\oV$-\icat{}. Then there is a natural
  equivalence of \icats{} between
  \begin{itemize}
  \item   lax symmetric monoidal $\oV$-functors $\phi_{*}\bar{\oC} \to
    \eD$,
  \item and oplax triangles
        \[
    \begin{tikzcd}
     \oC \ar[rr] \ar[dr, "\phi"', ""{name=A}] & & U\eD \ar[dl, "{\eD(\bbone,\blank)}"] \\
     & \oV
     \arrow[from=A,to=1-3,Rightarrow,shorten <=8pt, shorten >=20pt]
    \end{tikzcd}
  \]
  of lax symmetric monoidal functors.
  \end{itemize}
\end{corint}

\subsection*{Motivation}
\label{mot}

Our interest in \cref{thm:main V} mainly stems from its use in our
proof \cite{unfolding} of Lurie's ``unfolding'' equivalence \cite{LurieCob} that identifies symmetric
monoidal \incats{} with duals for objects and certain adjoints as
``chain complexes'' of symmetric monoidal \icats{} with duals: If $\eC$
is a symmetric monoidal \incat{} with duals, then it is in particular
closed, and so can be recovered via \cref{thm:main V} from the lax symmetric monoidal
functor
\[ \eC(\bbone, \blank) \colon \core[1] \eC \to \Cat_{(\infty,n-1)}\]
from its underlying \icat{}; we can unstraighten this to a fibration
\[ \mP\eC \to \core[1] \eC\] of symmetric monoidal
$(\infty,n-1)$-categories. This already gives the case $n=2$ of the
unfolding equivalence, and forms the first step in the general proof.

We briefly mention two other examples of enriched
\icats{} that can be constructed from \cref{thm:main V}:
\begin{itemize}
\item Let $\tMor(\oV)$ denote the Morita \itcat{} of a \SM{} \icat{}
  $\oV$ as defined in
  \cite{nmorita}; this has associative algebras in $\oV$ as objects,
  bimodules as morphisms (with composition given by relative tensor
  products), bimodule homomorphisms as 2-morphisms, and inherits a
  symmetric monoidal structure from that on $\oV$. The mapping
  functor
  \[ \tMor(\oV)^{\op} \times \tMor(\oV) \to \CATI\] is a lax symmetric
  monoidal functor of \itcats{}. Restricting this in the first
  variable to endomorphisms of the monoidal unit, we get a lax
  symmetric monoidal functor
  \[ \tMor(\oV) \to \FUN(B\oV, \CATI),\] where the target is the
  \itcat{} of $\oV$-modules in $\CATI$. This functor takes $A$ to
  $\Mod_{A}(\oV)$, and so takes values in $\oV$-modules in $\PrL$, which we
  can identify with presentable $\oV$-\icats{}. We thus in particular
  get a lax symmetric monoidal functor from $\core[1] \tMor(\oV)$ to
  $\oV$-\icats{}. This corresponds via \cref{thm:main V} to an
  enhancement of $\tMor(\oV)$ to a $\oV$-\itcat{}
  $\ecatname{M}{or}(\oV)$. We note that a version of the ``Deligne
  conjecture'' \cite{KontsevichSoibelman,McClureSmith} is an immediate consequence of this construction: the
  Hochschild cohomology of an associative algebra $A$ is the
  endomorphisms of $A$ as an $A$-module, which is the endomorphisms of
  the identity of $A$ viewed as an object in
  $\ecatname{M}{or}(\oV)$. It is therefore an $E_{2}$-algebra in $\oV$
  as the endomorphisms of an identity morphism is always an
  $E_{2}$-algebra for any $\oV$-\itcat{} \cite[Corollary 6.3.11]{enr}.
\item A \emph{6-functor formalism} in the sense of Mann--Scholze
  \cite{Scholze6} is a lax symmetric monoidal functor
  $D \colon \Span(\oC) \to \PrL$ from an \icat{} of spans (where $\oC$
  is typically some \icat{} of geometric objects, and there may be
  restrictions on the types of morphisms that appear in the
  spans). Here $\Span(\oC)$ inherits a symmetric monoidal structure
  from the cartesian product in $\oC$, for which every object is
  self-dual; thus $\Span(\oC)$ is closed with the internal Hom from
  $X$ to $Y$ given by $X \times Y$. The functor $D$ therefore
  corresponds via \cref{thm:main V} to a closed symmetric monoidal
  \itcat{} whose objects are those of $\oC$ and with mapping \icats{}
  from $X$ to $Y$ given by $D(X \times Y)$. This is the \emph{\itcat{}
    of kernels} studied by Heyer--Mann~\cite{HeyerMann}. The
  cocompletion of a particular \itcat{} of kernels for $\ell$-adic
  sheaves has also been studied by Gaitsgory, Rozenblyum, and
  Varshavsky~\cite{GRV} under the name \catname{AGCat}.
\end{itemize}

\subsection*{Overview}

In \cref{sec:enr weak alg} we introduce a new model for enriched
\iopds{}, which is an operadic analogue of Lurie's definition of
enriched \icats{} as weak modules. We then look at symmetric monoidal
enriched \icats{} within this model in \cref{sec:sm enr conj} and
further specialize to closed symmetric monoidal enrichments in
\cref{sec:closed}. In \cref{sec:lax tr} we then show that for a closed
symmetric monoidal $\oV$-\icat{} $\eC$ with underlying symmetric
monoidal \icat{} $\oC$, there is a canonical morphism
$\overline{\oC} \to \eC$ over a lax symmetric monoidal functor
$\eC(\bbone, \blank) \colon \oC \to \oV$ given by maps out of the
monoidal unit. The ability to construct these morphisms is the reason we
introduce our new model of enriched \iopds{}. We further show that
these maps give a lax natural transformation, which we then use in
\cref{sec:extend adjn} to construct an adjunction that restricts to
the equivalence of \cref{thm:main V}. To define this adjunction we use
a construction involving free fibrations that we discuss in
\cref{app:fibs}, after briefly reviewing fibrations of \itcats{}.

\subsection*{Notation}
We assume the reader is familiar with the basic concepts of \iopds{}
(in the sense of \cite{HA}) and \itcats{}.

\begin{itemize}
\item We denote generic \icats{} by $\oA, \oB, \oC, \dots$, \itcats{}
  by $\tA, \tB, \tC, \dots$ and
  enriched \icats{} by $\eA, \eB, \eC, \dots$.
\item $\SMCatl$ is the \icat{} of symmetric
  monoidal \icats{} and lax symmetric monoidal functors, and $\SMCATl$
  is the \itcat{} thereof.
\end{itemize}

\section{Enriched $\infty$-operads as weak algebras}
\label{sec:enr weak alg}

In this section we briefly introduce a new definition of enriched
\iopds{}. This is an operadic version of the approach to enriched
\icats{} originally due to Lurie~\cite[\S 4.7.1]{HA} and extensively developed by
Heine~\cite{HeineEnr,HeineBienr,HeineColim}; we start by briefly
recalling this definition:

If $\RM$ denotes the
\iopd{} for right modules over an associative algebra, an \iopd{}
$\oO \to \RM$ exhibits an \icat{} $\oM := \oO_{m}$ (where $m \in
\oO_{\angled{1}}$ is the object that gives the module) as a ``weak module'' over a
non-symmetric \iopd{} $\oP := a^{*}\oO$ (where $a \colon
\catname{Ass}\to \RM$ is the map from
the associative operad that extracts the algebra), in the sense that we
can talk about an \igpd{} $\oO_{\act}((m,p_{1},\dots,p_{n}),m')$ of
multimorphisms
\[ (m,p_{1},\dots,p_{n}) \to m'\]
for $p_{i} \in \oP_{\angled{1}}$ and $m,m' \in \oM$. Lurie and Heine define
enriched \icats{} as such \iopds{} over $\RM$ satisfying certain
properties:
\begin{itemize}
\item $\oP$ is a monoidal \icat{} $\oV^{\otimes}$,
\item for $v_{i} \in \oV$ and $m \in \oM$, the cocartesian map $(m,v_{1},\dots,v_{n}) \to (m,v_{1} \otimes \cdots \otimes v_{n})$ gives an equivalence
  \[ \oO_{\act}((m,v_{1} \otimes \cdots \otimes v_{n}),m') \isoto \oO_{\act}((m,v_{1},\dots,v_{n}),m')\]
  for any $m' \in \oM$,
\item the presheaf $v \mapsto \oO_{\act}((m,v),m')$ is representable
  by an object $\eM(m,m')$ of $\oV$ for all $m,m'$ in
  $\oM$ (which gives the $\oV$-enriched mapping object).
\end{itemize}
We can extract, for example, the enriched composition map
\[ \eM(m,m') \otimes \eM(m',m'') \to \eM(m,m'')\]
from the equivalences
\[
  \begin{split}
    \oV(\eM(m,m') \otimes \eM(m',m''), \eM(m,m'')) & \simeq
    \oO_{\act}((m,\eM(m,m') \otimes \eM(m',m'')), m'') \\ & \simeq 
    \oO_{\act}((m,\eM(m,m'), \eM(m',m'')), m''),
  \end{split}
\]
where on the right we have a canonical map, namely the composite
\[ (m, \eM(m,m'), \eM(m',m'')) \to (m', \eM(m',m'')) \to m''.\]

We will replace the base $\RM$ with an \iopd{} $\CA$ that encodes one
commutative algebra acting on another, and define
enriched \iopds{} as certain \iopds{} over $\CA$ satisfying analogues
of the conditions above.

\begin{defn}
  Let $\CA \to \xF_{*}$ be the \iopd{} that has
  \begin{itemize}
  \item as objects pairs $(\angled{n}, S \subseteq \angled{n}^{\circ})$,
  \item as morphisms $(\angled{n}, S) \to (\angled{m}, T)$ the
    morphisms $\phi \colon \angled{n} \to \angled{m}$ in $\xF_{*}$
    such that $\phi$ restricts to a map
    $S \cap \phi^{-1}(\angled{m}^{\circ}) \to T$.
  \end{itemize}
  The structure map to $\xF_{*}$ is the obvious projection.
\end{defn}

\begin{observation}
  $\CA$ is equivalent to $[1]^{\amalg}$, meaning the symmetric
  monoidal structure on $[1]$ given by the coproduct (= maximum):
  \begin{itemize}
  \item An object of $[1]^{\amalg}$ is a
    list $(\angled{n}, i_{1},\ldots,i_{n})$ where $i_{s} \in \{0,1\}$,
    which we identify with $(\angled{n}, S)$ with $S = \{s : i_{s} =
    1\}$.
  \item A morphism $(\angled{n},
    i_{1},\ldots,i_{n}) \to (\angled{m}, j_{1},\ldots,j_{m})$ in
    $[1]^{\amalg}$ is (since $[1]$ is a poset) given by a
    morphism $\phi \colon \angled{n} \to \angled{m}$ such that $i_{s}
    \leq j_{\phi(s)}$ if $\phi(s) \in \angled{m}^{\circ}$. If the
    corresponding subsets are $S \subseteq \angled{n}^{\circ}$ and $T
    \subseteq \angled{m}^{\circ}$, then this condition says precisely
    that $\phi(S \cap \phi^{-1}(\angled{m}^{\circ})) \subseteq T$.
  \end{itemize}
\end{observation}
We can therefore identify $\CA = [1]^{\amalg}$ as the Boardman--Vogt
tensor product $\xF_{*}\otimes [1]$ by \cite{HinichBVTens}*{Theorem
  2.3.1}. A $\CA$-algebra in a symmetric monoidal \icat{} $\oV$ is
hence equivalent to a morphism of commutative algebras $A \to B$
in $\oV$; we think of this as presenting an action of $A$ on $B$ as a
commutative algebra.

\begin{notation}
  We have two inclusions
  \[ e,u \colon \xF_{*} \to \CA,\] given by
  $\angled{n} \mapsto (\angled{n}, \emptyset)$ and
  $(\angled{n}, \angled{n}^{\circ})$. (For enriched \iopds{}, pulling
  back along these
  will extract, respectively, the symmetric monoidal \icat{}
  $e^{*}\uQ$ we enrich in and the underlying ordinary \iopd{}
  $u^{*}\uQ$ of the enriched \iopd{}.) 
\end{notation}

\begin{defn}
  An \emph{enriched \iopd{}} is an \iopd{} $\oP$ over $\CA$ such that
  the following conditions hold:
  \begin{enumerate}[(1)]
  \item $e^{*}\oP \to \xF_{*}$ is a symmetric monoidal \icat{}. We
    denote this as $\oV^{\otimes}$ from now on, and also abbreviate
    $\oO := u^{*}\oP$.
  \item For $v_{1},\dots,v_{n} \in \oV$ and $o_{1},\dots,o_{m} \in
    \oO_{\angled{1}}$, composition with the cocartesian morphism
    \[ (v_{1},\dots,v_{n},o_{1},\dots,o_{m}) \to (v_{1} \otimes \cdots
      \otimes v_{n}, o_{1},\dots,o_{m})\]
    gives an equivalence
    \[ \oP_{\act}((v_{1} \otimes \cdots
      \otimes v_{n}, o_{1},\dots,o_{m}), o') \isoto
      \oP_{\act}((v_{1},\dots,v_{n},o_{1},\dots,o_{m}), o'),\]
    for all $o' \in \oO_{\angled{1}}$.
  \item For all $o_{1},\ldots,o_{n},o' \in \uO_{\angled{1}}$, the
    presheaf $v \mapsto \oP_{\act}((v, o_{1},\ldots,o_{n}), o')$ is
    representable by an object $\eO(o_{1},\ldots,o_{n}; o') \in \uV$.
  \end{enumerate}
  We define the \itcat{} $\ENROPD$ as the full
  subcategory of $\OPDsl{\CA}$ spanned by the enriched \iopds{},
  where $\OPD$ is the \itcat{} of \iopds{}. The
  functor $e^{*}$ of pullback along $e$ then gives a functor
  \[ p \colon \ENROPD \to \SMCATl,\] where
  $\SMCATl \subseteq \OPD$ is the full sub-\itcat{} of symmetric
  monoidal \icats{} and lax symmetric monoidal functors. For $\oV$ a
  symmetric monoidal \icat{}, the \itcat{} $\ENROPD(\oV)$ of
  $\oV$-enriched \iopds{} is the fibre of $p$ at $\oV$.
\end{defn}

\begin{notation}
  We will often denote objects of $\ENROPD$ as $\uract{O}{V}$ where
  $\oO = u^{*}\uract{O}{V}$ and
  $\oV^{\otimes} = e^{*}\uract{O}{V}$, when the particular definition
  of enriched \iopds{} above is relevant. Similarly, we will denote
  morphisms in $\ENROPD$ as
  $\ract{F}{\phi} \colon \uract{O}{V} \to \uract{P}{W}$, where
  $F = u^{*}\ract{F}{\phi}$ and $\phi = e^{*}\ract{F}{\phi}$.
\end{notation}

\begin{observation}
  Suppose $\uract{O}{V}$ is an enriched \iopd{}. Then the enriched
  mapping objects for $o_{1},\dots,o_{n},o'$ in $\oO_{\angled{1}}$
  satisfy
  \[
    \begin{split}
      \oV(\bbone, \eO(o_{1},\dots,o_{n};o')) & \simeq
                                               (\uract{O}{V})^{\act}((\bbone, o_{1},\dots,o_{n}), o') \\
      & \simeq \oO^{\act}((o_{1},\dots,o_{n}), o').
    \end{split}
  \]
\end{observation}

We assume that our enriched \iopds{} are 2-functorial under change of
enrichment, in the following sense:
\begin{assumption}\label{assumption:eopdfib}
  The functor $p \colon \ENROPD \to \SMCATl$ is a $(0,1)$-fibration of
  \itcats{} (see \cref{def:01fib}). A morphism $\ract{F}{\phi} \colon \uract{O}{V} \to
  \uract{P}{W}$ is $p$-cocartesian \IFF{}
  \begin{itemize}
  \item $F_{\angled{1}} \colon \oO_{\angled{1}} \to \oP_{\angled{1}}$ is essentially surjective,
  \item for all $o_{1},\dots,o_{n},o' \oO_{\angled{1}}$, the canonical
    map\footnote{$\ract{F}{\phi}$ gives a map 
      $ (\uract{O}{V})^{\act}((v, o_{1},\dots,o_{n}), o') \to
        (\uract{P}{W})^{\act}((\phi v, F o_{1},\dots, F o_{n}), F o')
      $
      where the source is equivalent to
      $\oV(v, \eO(o_{1},\dots,o_{n};o'))$ and the target to
      $\oW(\phi v, \eP(F o_{1},\dots,F o_{n}; F o'))$. For
      $v = \eO(o_{1},\dots,o_{n};o')$, the identity gives our
      canonical map.}
    \[ \phi \eO(o_{1},\dots,o_{n};o')  \to \eP(F o_{1},\dots,F o_{n}; F
      o')\]
    is an equivalence.
  \end{itemize}
  A 2-morphism $\alpha \colon \ract{F}{\phi} \to \ract{G}{\psi}$ is
  $p$-cartesian \IFF{} $u^{*}\alpha \colon F \to G$ is an equivalence.
\end{assumption}

\begin{remark}
  Under the straightening equivalence for $(0,1)$-fibrations,
  \cref{assumption:eopdfib} says that the assignment $\oV
  \mapsto \ENROPD(\oV)$ determines a functor of \itcats{} $\SMCATl \to
  \CATIT$. For $\phi \colon \oV \to \oW$, the characterization of
  cocartesian morphisms above says that the pushforward of $\eO \in
  \ENROPD(\oV)$ is a $\oW$-\iopd{} $\phi_{*}\eO$ with the same objects
  as $\eO$ and mapping objects $\phi \eO(\blank;\blank)$, as we would
  expect. For a natural transformation $\alpha \colon \phi \to \psi$
  we have a natural transformation $\alpha_{*} \colon \phi_{*} \to \psi_{*}$ given by
  applying $\alpha$ on mapping objects. To relate this to the
  characterization of cartesian 2-morphisms, recall that a 2-morphism
  between 1-morphisms $\eO \to \eP$ over $\alpha$ corresponds to a triangle of $\oW$-\iopds{}
  \[
    \begin{tikzcd}[row sep=tiny]
      \psi_{*}\eO \ar[dr, ""{name=A}'] \ar[dd, "\alpha_{*,\eO}"'] \\
      & \eP \\
      \phi_{*}\eO; \ar[ur]
      \ar[from=A,to=3-1,Rightarrow,shorten <= 8pt, shorten >= 8pt
      ]
    \end{tikzcd}
  \]
  this 2-morphism is cartesian \IFF{} the triangle commutes. Thus the
  claim amounts to the observation that if a 2-morphism of \iopds{}
  over $\CA$ becomes invertible after applying both $e^{*}$ and
  $u^{*}$, then it is invertible.
\end{remark}

\begin{remark}\label{rmk:assumption sketch}
  While we will not go into the details here, we give a brief sketch
  of how one can prove \cref{assumption:eopdfib}:
  \begin{itemize}
  \item To construct a $p$-cocartesian morphism in $\ENROPD$ from
    $\uract{O}{V}$ over
    $\phi \colon \oV \to \oW$, we can first take the envelope of
    $\uract{O}{V}$ over $\CA$ to obtain a $\CA$-algebra in $\CatI$ and
    then take presheaves to get a $\CA$-algebra in $\PrL$, \ie{} an
    algebra in $\PrL$ for presheaves on $\Env(\oV)$.\footnote{Note
      that if we also remember the objects that come from $\oO$, this
      produces an \iopd{} enriched in $\PSh(\Env(\oV))$ in
      the sense of Zetto~\cite{ZettoOpd}.}
    We can then take the relative tensor product with
    $\PSh(\Env(\oW))$ viewed as a commutative algebra over
    $\PSh(\Env(\oV))$ via $\PSh(\Env(\phi))$ and restrict to the image
    of the objects from $\oO$ to obtain the $p$-cocartesian
    morphism. This is an operadic analogue of the construction of
    change of enrichment for enriched \icats{} in
    \cite{HeineBienr,ReutterZetto}.
  \item To construct $p$-cartesian 2-morphisms, we observe that for
    any \iopd{} $\oO$ over $\CA$, the functor $e^{*}\oO \to \oO$ has a
    left adjoint $\ell \colon \oO \to e^{*}\oO$. 
    Moreover, $\oO$ has
    $\ell$-cartesian lifts of all morphisms that map to equivalences
    in $\xF_{*}$, and these are precisely the morphisms that map to
    equivalences in $\CA$ and whose components in
    $u^{*}\oO_{\angled{1}}$ are all invertible. This allows us to
    construct $\ell$-cartesian lifts of natural transformations, which
    give $e^{*}$-cartesian 2-morphisms.
  \end{itemize}
\end{remark}

\section{Symmetric monoidal enriched
  $\infty$-categories}
\label{sec:sm enr conj}

We now define symmetric monoidal enriched \icats{} as a full
subcategory of enriched \iopds{}.

\begin{defn}
  An enriched \iopd{} $\uract{O}{V}$ is a \emph{symmetric monoidal
    enriched \icat{}} if
  \begin{enumerate}[(1)]
  \item its underlying \iopd{} $\oO$ is a symmetric monoidal
    \icat{},
  \item for all
    $v \in \oV, x_{1},\ldots,x_{n},y \in \oO_{\angled{1}}$, the
    cocartesian morphism $(x_{1},\dots,x_{n}) \to x_{1} \otimes \cdots
    \otimes x_{n}$ induces an equivalence
    \[ (\uract{O}{V})_{\act}((v, x_{1} \otimes \cdots \otimes x_{n}),
      y) \isoto (\uract{O}{V})_{\act}((v, x_{1}, \dots, x_{n}),
      y);\]
    equivalently, the morphism of enriched mapping objects
    \[ \eO(x_{1}\otimes \cdots \otimes x_{n}; y) \isoto
      \eO(x_{1},\ldots,x_{n}; y)\]
    is an equivalence.
  \end{enumerate}
  We write $\SMENRl$ for the full sub-\itcat{} of $\ENROPD$ spanned by
  the symmetric monoidal enriched \icats{},
  $p' \colon \SMENRl \to \SMCATl$ for the restriction of $p = e^{*}$,
  and $\SMENRl(\oV)$ for the fibre of $p'$ at a symmetric monoidal
  \icat{} $\oV$.
\end{defn}

\begin{observation}\label{obs:esm simplify}
  The preceding definition is slightly redundant: Suppose
  $\uract{O}{V}$ is an enriched \iopd{} so that for all
  $x_{1},\dots,x_{n}$ in $\oO_{\angled{1}}$ there exists an active
  morphism $(x_{1},\dots,x_{n}) \to x'$ with the property that for all
  $y \in \oO$ the induced map
  \[ \eO(x'; y) \to
    \eO(x_{1},\ldots,x_{n}; y)\]
  is an equivalence, or equivalently the maps
  \[ (\uract{O}{V})_{\act}((v, x'),
    y) \isoto (\uract{O}{V})_{\act}((v, x_{1}, \dots, x_{n}),
    y)\]
  are equivalences for all $v \in \oV$. Then $\oO$ is a symmetric
  monoidal \icat{}, since taking $v = \bbone$ we get the required
  equivalence
  \[ \oO_{\act}(x', y) \isoto\oO_{\act}((x_{1},\dots,x_{n}), y).\]
\end{observation}

\begin{lemma}\label{lem:SME cond}
  An \iopd{} $\oP$ over $\CA$ is a symmetric monoidal enriched \icat{}
  \IFF{} the following conditions hold:
  \begin{enumerate}[(1)]
  \item the \iopds{} $\oV^{\otimes} = e^{*}\oP$ and $\oC^{\otimes} = u^{*}\oP$ are symmetric monoidal
    \icats{},
  \item for $v_{1},\dots,v_{n} \in \oV$ and $c_{1},\dots,c_{m} \in
    \oC$, the cocartesian morphisms $(v_{1},\dots,v_{n}) \to v_{1}
    \otimes \cdots v_{n}$ and $(c_{1},\dots,c_{m}) \to c_{1} \otimes
    \cdots \otimes c_{m}$ induce for all $c' \in \oC$ an equivalence
    \[ \oP_{\act}((v_{1} \otimes \cdots \otimes v_{n}, c_{1} \otimes
      \cdots \otimes c_{m}),c') \isoto
      \oP_{\act}((v_{1},\dots,v_{n},c_{1},\dots,c_{m}),c').\]
  \item The presheaf $v \mapsto \oP_{\act}((v,c), c')$ is
    representable for all $c,c' \in \oC$.
  \end{enumerate}
\end{lemma}
\begin{proof}
  If the given conditions hold, it suffices to prove that $\oP$ is an
  enriched \iopd{}. In (2), note that we have a commutative square
  \[
    \begin{tikzcd}
     \oP_{\act}((v_{1} \otimes \cdots \otimes v_{n}, c_{1} \otimes
      \cdots \otimes c_{m}),c') \ar[r]
     \ar[d] & \oP_{\act}((v_{1} \otimes \cdots \otimes v_{n},c_{1},\dots,c_{m}),c')  \ar[d] \\
     \oP_{\act}((v_{1},\dots,v_{n},c_{1} \otimes \cdots \otimes c_{m}),c') \ar[r] & \oP_{\act}((v_{1},\dots,v_{n},c_{1},\dots,c_{m}),c').
    \end{tikzcd}
  \]
  Here (2) implies that
  the top horizontal and left vertical
  morphisms, as well as the composite in the square, are equivalences,
  so the right vertical and bottom horizontal morphisms are also
  equivalences. This gives the second condition in the definition of
  an enriched \iopd{}, while the presheaf
  $ \oP_{\act}((\blank, c_{1},\dots,c_{m}), c')$ is equivalent to
  $\oP_{\act}((\blank, c_{1}\otimes \cdots \otimes c_{m}), c')$, and so
  is representable as the latter is.

  Conversely, if $\oP$ is a symmetric monoidal enriched \icat{}, we
  need to show that condition (2) here holds. In the square above, the
  vertical morphisms are equivalences since $\oP$ is an enriched
  \iopd{}, while the top horizontal morphisms is an equivalence since
  it is moreover a symmetric monoidal enriched \icat{}. The bottom
  horizontal morphism and the composite in the square are therefore
  also equivalences.
\end{proof}

\begin{propn}\label{propn:SMENRl cocart}
  The functor $p' \colon \SMENRl \to \SMCATl$ is a
  $(0,1)$-fibration. A morphism $\ract{F}{\phi} \colon \utract{C}{V}
  \to \utract{D}{V}$ is $p'$-cocartesian \IFF{}
  \begin{itemize}
  \item $F$ is essentially surjective and strong symmetric monoidal,
  \item for all $c,c' \in \oC$, the canonical
    map
    \[ \phi \eC(c,c') \to \eD(F c, F c'),\]
    induced by taking $v = \eC(c,c')$ in
    \[ (\utract{C}{V})^{\act}((v, c), c') \to
      (\utract{D}{W})^{\act}((\phi v, F c), F c'),\]
    is an equivalence.
  \end{itemize}
  A 2-morphism if $p'$-cartesian \IFF{} its image under $u^{*}$ is an
  equivalence.
\end{propn}
\begin{proof}
  Since $p$ is a $(0,1)$-fibration by \cref{assumption:eopdfib}, to
  see that the same holds for its restriction $p'$ it suffices to
  check that if $\utract{C}{V}$ is a symmetric monoidal enriched
  \icat{} and $\utract{C}{V} \to \phi_{*}\utract{C}{V}$ is the
  $p$-cocartesian morphism over $\phi \colon \oV \to \oW$, then
  $\uract{O}{W} := \phi_{*}\utract{C}{V}$ is again a symmetric
  monoidal enriched \icat{}.

  Given $x_{1},\dots,x_{n},x'$ in $\oO_{\angled{1}}$, since $F$ is
  essentially surjective we can find $c_{i},c' \in \oC$ with $F(c_{i})
  \simeq x_{i}$, $F(c')\simeq x'$. Then the cocartesian morphism
  $(c_{1}, \dots, c_{n}) \to c_{1} \otimes \dots \otimes c_{n}$
  gives a commutative square
  \[
    \begin{tikzcd}
     \phi \eC(c_{1}\otimes \cdots \otimes c_{n}, c') \ar[r, "\sim"]
     \ar[d, "\sim"'] & \eO(F(c_{1} \otimes \cdots \otimes c_{n}), y)  \ar[d] \\
     \phi \eC^{\otimes}(c_{1},\dots,c_{n};c') \ar[r, "\sim"'] &
     \eO(x_{1},\dots,x_{n}; y).
    \end{tikzcd}
  \]
  This implies that the morphism
  $(x_{1},\dots,x_{n}) \to F(c_{1} \otimes \cdots \otimes c_{n})$ has the required
  property from \cref{obs:esm simplify} for $\uract{O}{W}$ to be a
  symmetric monoidal enriched \icat{}. Moreover, $F$ is strong
  symmetric monoidal, so the stated conditions for a $p'$-cocartesian
  morphism are necessary.

  To see that the conditions for $\ract{F}{\phi} \colon \utract{C}{V}
  \to \utract{D}{V}$ to be $p'$-cocartesian are sufficient, observe
  that for $c_{1},\dots,c_{n},c'$ in $\oC$ we have a commutative square
  \[
    \begin{tikzcd}
     \phi \eC(c_{1} \otimes \cdots \otimes c_{n}, c') \ar[r]
     \ar[d] & \eD(F(c_{1} \otimes \cdots \otimes c_{n}), Fc')  \ar[d] \\
     \phi \eC^{\otimes}(c_{1},\dots,c_{n}; c') \ar[r] &
     \eD^{\otimes}(Fc_{1},\dots,Fc_{n}; Fc').
    \end{tikzcd}
  \]
  Here the top horizontal morphism is an equivalence by assumption,
  and the left vertical morphism is an equivalence since
  $\utract{C}{V}$ is a symmetric monoidal enriched \icat{}. The right
  vertical morphism is also an equivalence since $F$ is by assumption
  strong symmetric monoidal, so the bottom horizontal morphism must
  also be an equivalence, as required.
\end{proof}

\section{Closed symmetric monoidal structures}\label{sec:closed}

We now define closed symmetric monoidal enriched \icats{} and consider
the self-enrichment of closed symmetric monoidal \icats{}.

\begin{defn}
  We say that a symmetric monoidal enriched \icat{} $\utract{C}{V}$ is
  \emph{closed} if for all objects $x,y \in \oC$ there exists an
  object $[x,y]$ and a morphism $\ev \colon [x,y] \otimes x \to y$ such that
  \[ \eC(z, [x,y]) \to \eC(z \otimes x, [x,y] \otimes x) \xto{ev_{*}} \eC(z
    \otimes x, y)\]
  is an equivalence for all $z \in \oC$. We write $\SMENRlcl$ for the
  full sub-\itcat{} of $\SMENRl$ spanned by the closed symmetric
  monoidal enriched \icats{}, and
  \[ P \colon \SMENRlcl \to \SMCATl\]
  for the restriction of $p'$ to this full sub-\itcat{}.
\end{defn}

\begin{notation}
  We write $U \colon \SMENRlcl \to \SMCATlcl$ for the restriction of
  $u^{*} \colon \SMENRl \to \SMCATl$.
\end{notation}

\begin{lemma}
  A symmetric monoidal enriched \icat{} $\utract{C}{V}$ is closed
  \IFF{} the underlying symmetric monoidal \icat{} $\oC$ is closed,
  and for all $x,y \in \oC$ the map
  \[ \eC(\bbone, [x,y]) \to \eC(x, [x,y] \otimes x) \xto{\ev_{*}} \eC(x, y)\]
  is an equivalence.
\end{lemma}
\begin{proof}
  We only need to check that the given condition is sufficient, which
  follows by unpacking the composite equivalence
  \[ \eC(z, [x,y]) \simeq \eC(\bbone, [z, [x,y]]) \simeq \eC(\bbone, [z
    \otimes x, y]) \simeq \eC(z \otimes x, y)\] to see it is given by
  tensoring with $x$ and composing with evaluation.
\end{proof}

\begin{propn}\label{propn:SMENRlcl cocart}
  The functor $P \colon \SMENRlcl \to \SMCATl$ is a
  $(0,1)$-fibration. A morphism
  $\ract{F}{\phi} \colon \utract{C}{V} \to \utract{D}{W}$ is
  $P$-cocartesian \IFF{}
  \begin{itemize}
  \item   $F$ is essentially surjective and strong
    symmetric monoidal,
  \item $F$ preserves internal Homs, \ie{} the canonical map
    \[ F[x,y]_{\oC} \to [Fx,Fy]_{\oD}\]
    is an equivalence,
  \item for every
  object $c \in \oC$, the morphism
  \[ \phi\eC(\bbone_{\oC}, c) \to \eD(\bbone_{\oD}, F(c))\] is an
  equivalence in $\oW$.
  \end{itemize}
  A 2-morphism is $P$-cartesian \IFF{} its image
  under $U$ is an equivalence.
\end{propn}
\begin{proof}
  Since $P$ is the restriction of the $(0,1)$-fibration $p'$ to a full
  sub-\itcat{}, to see that $P$ is a $(0,1)$-fibration it suffices to check that if
  $\utract{C}{V}$ is closed, then its cocartesian transport
  $\phi_{*}(\utract{C}{V})$ is again
  closed for any lax symmetric monoidal functor
  $\phi \colon \uV \to \uW$. More precisely,
  the image of the internal Hom in $\oC$ remains one in
  $\phi_{*}(\utract{C}{V})$, which is clear from the commutative
  diagram
  \[
    \begin{tikzcd}
      \phi\eC(z, [x,y]) \ar[r, "\sim"]  \ar[ddr, "\sim"', bend right] & \phi\eC(z, [x,y]) \otimes \bbone_{\uW} \ar[d] \ar[dr] \\
      & \phi\eC(z, [x,y]) \otimes \phi \bbone_{\uV} \ar[r]  \ar[d] & \phi \eC(z, [x,y]) \otimes \phi \eC(x,x) \ar[d] \\
      & \phi (\eC(z, [x,y]) \otimes \bbone_{\uV}) \ar[r] & \phi (\eC(z, [x,y]) \otimes \eC(x,x)) \ar[d] \\
      & & \phi \eC(z \otimes x, x \otimes [x,y]) \ar[d] \\
       & & \phi \eC(z \otimes x, y),
    \end{tikzcd}
  \]
  where the composite along the top is the one we must check is an
  equivalence, and the composite along the bottom is $\phi$ applied to
  the equivalence witnessing that $\utract{C}{V}$ is
  $\uV$-closed. This also shows that a cocartesian morphism
  necessarily preserves internal Homs, so the given conditions are
  necessary.

  It remains to show that the claimed characterization of
  $P$-cocartesian
  morphisms implies the conditions for $p'$-cocartesian morphisms in \cref{propn:SMENRl
    cocart}. Assuming that $F$ is strong symmetric monoidal and preserves internal Homs, we have for
  $x, y \in \oC$ a
  commutative square
  \[
    \begin{tikzcd}
     \phi \eC(\bbone, [x,y]) \ar[r] \ar[d, "\sim"'] & \eD(\bbone, F[x,y])  \ar[d, "\sim"] \\
     \phi \eC(x,y) \ar[r] & \eD(Fx,Fy)
    \end{tikzcd}
  \]
  so the bottom morphism is an equivalence if the top one is
  invertible. This shows that the conditions are sufficient.
\end{proof}

\begin{propn}
  Let $\pi^{*} \colon \OPD \to \OPDsl{\CA}$ be the functor given by
  pullback along the structure map $\pi \colon \CA \to \xF_{*}$. Then
  $\pi^{*}\oO$ is an enriched \iopd{} \IFF{} the \iopd{} $\oO$ is a
  closed symmetric monoidal \icat{}, in which case $\pi^{*}\oO$ is a
  closed symmetric monoidal enriched \icat{}.
\end{propn}
\begin{proof}
  Since $e^{*}\pi^{*}\oO \simeq \oO$, we clearly need $\oO$ to be
  symmetric monoidal. The second condition in \cref{lem:SME cond} is
  then automatic, as $\pi^{*}\oO \to \CA$ is a cocartesian fibration,
  while the third condition asks for the presheaf
  \[ z \mapsto \oC(z \otimes x, y)\] to be representable for all
  $x,y$, \ie{} for $\oO_{\angled{1}}$ to be closed. The condition for
  $\pi^{*}\oO$ to be closed then asks for
  \[ [z, [x,y]] \to [z \otimes x, [x,y] \otimes x] \to [z \otimes x,
    y]\] to be an equivalence, which follows from $\oO_{\angled{1}}$
  being closed.
\end{proof}

\begin{notation}
  Let $\SMCATlcl$ denote the full sub-\itcat{} of $\SMCATl$ spanned by
  the closed symmetric monoidal \icats{}. We write
  \[ S \colon \SMCATlcl \to \SMENRlcl \]
  for the functor given by $\pi^{*}$ on this full sub-\itcat{}. We
  often write $\bar{\oC}$ for $S(\oC^{\otimes})$; this is the
  \emph{self-enrichment} of the closed symmetric monoidal \icat{}
  $\oC$.
\end{notation}

\begin{observation}
  We have $US \simeq \id_{\SMCATlcl}$, as this is the restriction of
  $u^{*}\pi^{*} \simeq \id_{\OPD}$.
\end{observation}

\section{A lax transformation}\label{sec:lax tr}

Our goal in this section is to construct a lax natural transformation
$\alpha \colon SU \to \id_{\SMENRlcl}$. This will essentially
provide the counit of our desired adjunction in the next
section. Given a closed symmetric monoidal $\oV$-\icat{} $\eC$, the
component $SU(\eC) \to \eC$ will live over a lax symmetric monoidal
functor $U(\nC) \to \oV$ given by $\eC(\bbone,\blank)$. 

\begin{construction}
  Define a natural transformation $\tau \colon \CA \times [1] \to \CA$ by
  \[ \tau(\angled{n}, S, i) =
    \begin{cases}
      (\angled{n}, S), & i = 0 \\
      (\angled{n}, \angled{n}^{\circ}), & i = 1;
    \end{cases}
  \]
  thus $\tau_{0} = \id_{\CA}$ and $\tau_{1}$ is the composite
  \[ \CA \xto{\pi} \xF_{*} \xto{u} \CA \] Given
  $\utract{M}{V} \to \CA$ we can pull back along $\tau$ to get
  $\tau^{*}\utract{M}{V} \to \CA \times [1]$ with fibre
  $\utract{M}{V}$ at $0$ and $\utract{M}{M} := \pi^{*}\uM^{\otimes}$
  at $1$.
\end{construction}

\begin{propn}\label{extract}
  Suppose $\utract{M}{V}$ is a symmetric monoidal $\uV$-\icat{}. Then in the triangle
  \[
    \begin{tikzcd}
      \tau^{*}\utract{M}{V} \ar[rr] \ar[dr] & & \CA \times [1] \ar[dl] \\
      & {[1]},
    \end{tikzcd}
  \]
  the functor to $[1]$ is a cartesian fibration and the top horizontal
  functor preserves cartesian morphisms. 
\end{propn}
\begin{proof}
  Given an object $M = (m_{1},\ldots,m_{n})$ in
  \[ \tau^{*}\utract{M}{V}_{(\angled{n}, S, 1)} \simeq \uM^{\otimes}_{\angled{n}} \simeq \uM^{\times n},\]
  we must exhibit a cartesian morphism $\gamma_{M} \colon M' \to M$ to
  this object, which should lie over the identity of $(\angled{n},
  S)$. We define $M' = (m'_{1},\ldots,m'_{n})$  and $\gamma_{M}$ by
  \[ m'_{i} :=
    \begin{cases}
      m_{i}, & i \in S, \\
      \nM(\bbone_{\uM}, m_{i}), & i \notin S,
    \end{cases}
    \qquad\qquad (\gamma_{M})_{i} =
    \begin{cases}
      \id_{m_{i}}, & i \in S, \\
      q_{m_{i}}, & i \notin S,
    \end{cases}
  \]
  where  $q_{m} \colon \nM(\bbone_{\uM}, m) \to m$ for $m \in
  \oM$ is defined to be the map that corresponds to the identity of $\nM(\bbone_{\uM}, m)$ under the equivalences
  \[
    (\utract{M}{V})_{\act}(v, m) \simeq (\utract{M}{V})_{\act}((v, \bbone_{\uM}), m) \simeq \oV(v, \nM(\bbone_{\uM}, m)),
  \]
  for $v \in \uV$. Now we must show that for an object $X = (x_{1},\ldots,x_{m})$ over $(\angled{m}, T, 0)$, composition with $\gamma_{M}$ induces an equivalence
  \[ \gamma_{M,*} \colon \utract{M}{V}(X, M') \isoto \utract{M}{V}(X, M).\]
  This map lives over $\CA((\angled{m}, T), (\angled{n}, S))$, and it suffices it gives an equivalence on the fibres over each map $\phi \colon (\angled{m}, T) \to (\angled{n}, S)$ in $\CA$. Suppose
  \[ (\angled{m}, T) \xto{\alpha} (\angled{k}, U) \xto{\beta} (\angled{n}, S)\]
  is the inert--active factorization of $\phi$; if $\bar{\alpha} \colon X \to \alpha_{!}X$ is cocartesian over the inert morphism $\alpha$, then composition with $\bar{\alpha}$ gives a commutative square
  \[
    \begin{tikzcd}
      \utract{M}{V}(\alpha_{!}X, M')_{\beta} \ar[r, "\gamma_{M,*}"] \ar[d, "\bar{\alpha}^{*}"]  & \utract{M}{V}(X, M)_{\beta} \ar[d, "\bar{\alpha}^{*}"] \\
      \utract{M}{V}(X, M')_{\phi} \ar[r, "\gamma_{M,*}"] \ar[r] &  \utract{M}{V}(X, M)_{\phi},
    \end{tikzcd}
  \]
  where the vertical morphisms are equivalences. It therefore suffices to consider the case where $\phi$ is active. For $i = 1,\ldots,m$, let $\phi_{i} \colon (\angled{m_{i}}, T_{i}) \actto (\angled{1}, S_{i})$ be defined as the inert--active factorization in the square
  \[
    \begin{tikzcd}
      (\angled{m}, T) \ar[r, active] \ar[d, inert, "p_{i}"]& (\angled{n}, S) \ar[d, inert, "\rho_{i}"] \\
      (\angled{m_{i}}, T_{i}) \ar[r, active, "\phi_{i}"] & (\angled{1}, S_{i})
    \end{tikzcd}
  \]
  (where $S_{i} = \{i\} \cap S$). Again using the inert-active factorization system and the decomposition of mapping spaces in \iopds{} as products, we can identify the map we are interested in as
  \[ \prod_{i = 1}^{n}(\utract{M}{V})_{\act}(X_{i}, M'_{i}) \to \prod_{i = 1}^{n} (\utract{M}{V})_{\act}(X_{i}, M_{i}),\]
  where $X \to X_{i}$ is cocartesian over $p_{i}$. It thus suffices to consider the case where $n = 1$. If $M = (m)$ lives over $(\angled{1}, \{1\})$ then $\gamma_{M} = \id_{m}$ and there is nothing to prove. We are left with the case where $M$ lives over $(\angled{1}, \emptyset)$ and $\gamma_{M}=q_{m}$. Here we must have $T = \emptyset$ for a map $(\angled{m}, T ) \to (\angled{1}, \emptyset)$ to exist, so that $X = (v_{1},\ldots,v_{m})$ for $v_{i} \in \uV$. Then for $v = v_{1} \otimes \cdots \otimes v_{m}$ we have
  \[
    \begin{tikzcd}
      (\utract{M}{V})_{\act}(X, \nM(\bbone_{\uM}, m)) \ar[r, "q_{m,*}"] \ar[d, "\sim"] & (\utract{M}{V})_{\act}(X, m) \ar[d, "\sim"] \\
      \oV(v, \nM(\bbone_{\uM}, m)) \ar[r, "q_{m,*}"] & (\utract{M}{V})_{\act}(v, m),
    \end{tikzcd}
  \]
  where the bottom horizontal map is an equivalence by the definition of $q_{m}$.
\end{proof}

\begin{cor}
  Given a symmetric monoidal enriched \icat{} $\utract{M}{V}$ we can
  unstraighten $\tau^{*}\utract{M}{V}$ over $[1]$ to get a commutative triangle
\[
  \begin{tikzcd}
    \utract{M}{M} \ar[rr] \ar[dr] & & \utract{M}{V} \ar[dl] \\
    & \CA,
  \end{tikzcd}
\]
where the top horizontal functor preserves inert morphisms. In
particular, there is a canonical lax symmetric monoidal functor
$\uM \to \uV$ given by $\nM(\bbone_{\uM}, \blank)$.
\end{cor}
\begin{proof}
  It remains to show that the straightening preserves inert morphisms, but this amounts to checking that for $M \intto M_{0}$ cocartesian over an inert map $\alpha$ in $\CA$, we have a commutative square
  \[
    \begin{tikzcd}
      M' \ar[r, "q_{M}"] \ar[d, inert] & M \ar[d, inert] \\
      M'_{0} \ar[r, "q_{M_{0}}"] & M_{0},
    \end{tikzcd}
  \]
  where $M' \to M'_{0}$ is also cocartesian over $\alpha$, but this is
  clear from the definition since we are just projecting off some of
  the factors in $q_{M}$.
\end{proof}

\begin{observation}\label{obs:taucartmor}
  Note that while $\tau^{*}\utract{M}{V} \to [1]$ is a cartesian
  fibration, the functor
  $\tau^{*}\utract{M}{V} \to \tau^{*}\utract{N}{W}$ over $[1]$ induced
  by a morphism
  $\ract{F}{\phi} \colon \utract{M}{V} \to \utract{N}{W}$ does not in
  general preserve cartesian morphisms --- indeed, from our
  description of the latter we see that they are preserved \IFF{} the
  natural map
\[ \phi (\nM(\bbone_{\uM}, m)) \to \nN(\bbone_{\uN}, F(m)) \] is an
equivalence for all $m \in \uM$.
\end{observation}

\begin{construction}
  Arbitrary morphisms between cartesian fibrations correspond to oplax
  transformations by \cite{HHLN1}*{Theorem E}, unstraightening over
  $[1]$ and identifying $[1]^{\op}$ with $[1]$ gives after restricting
  to the closed case a functor of \itcats{}
  \begin{equation}
    \label{eq:SUlax}
    \alpha \colon \SMENRlcl \to \ARopl(\SMENRlcl),  
  \end{equation}
  given by $SU$ at $0$ and $\id$ at $1$,
  which we can also view as a lax natural transformation from $SU$ to
  $\id$.
\end{construction}

\begin{propn}\label{propn:alpha props}
  The lax transformation $\alpha$ has the following properties:
  \begin{enumerate}[(i)]
  \item The composite $U \simeq USU \xto{U\alpha} U$ is equivalent to
    $\id_{U}$.    
  \item The composite $S \simeq SUS \xto{\alpha S} S$ is equivalent to
    $\id_{S}$.
  \item The lax naturality squares of $\alpha$ all contain
    $P$-cartesian 2-morphisms.
  \item The lax naturality squares of $\alpha$ commute at
    $P$-cocartesian morphisms.
  \item $\alpha_{\utract{C}{V}}$ is a $P$-cocartesian morphism for all
    objects $\utract{C}{V} \in \SMENRlcl$.
  \end{enumerate}
\end{propn}
\begin{proof}
  By definition, $\alpha_{\utract{C}{V}}$ is the straightening of
  $\tau^{*}\utract{C}{V}$ over $[1]$. From the definition of $\tau$ we
  see immediately that the composite $\tau \circ (u \times \id)$ is
  the identity transformation on $u$, while $\pi \circ
  \tau$ is the identity transformation on $\pi$. Thus $(u \times
  \id)^{*}\tau^{*}\utract{C}{V}$ straightens to
  $\id_{U\utract{C}{V}}$ and $\tau^{*}\pi^{*}\oC^{\otimes}$
  straightens to the identity on $\pi^{*}\oC^{\otimes}$, which proves
  (i) and (ii).

  Since $P$-cartesian 2-morphisms are characterized as
  those that go to equivalences under $U$, (i) immediately implies
  (iii). To prove (iv) it then suffices to show that the image
  of the naturality 2-morphism for $\ract{F}{\phi} \colon
  \utract{C}{V} \to \utract{D}{W}$  goes to an equivalence under $P$
  when $\ract{F}{\phi}$ is a $P$-cocartesian morphism.   
  Unpacking the lax naturality square of $\alpha$, we see
  that this gives the natural map
  \[ \phi \eC(\bbone, \blank) \to \eD(\bbone, F(\blank)),\]
  which is indeed an equivalence for a $P$-cocartesian morphism by
  \cref{propn:SMENRlcl cocart}.
  
  It remains to prove (v). The first two conditions in
  \cref{propn:SMENRlcl cocart} are clear, since
  $U\alpha_{\utract{C}{V}}$ is an equivalence by (i). The third
  condition unpacks to
  \[ \eC(\bbone,[\bbone,c]) \to \eC(\bbone, c)\]
  being an equivalence, which is clear since $\utract{C}{V}$ is by
  assumption closed.
\end{proof}

\section{Extending to an adjunction}\label{sec:extend adjn}
Our goal in this section is to apply
\cref{propn:extend adjn} to extend the self-enrichment functor $S
\colon \SMCATlcl \to \SMENRlc$ to a left adjoint; this will then
quickly lead to the proof of \cref{thm:main V}. We have a commutative triangle of \itcats{}
\[
  \begin{tikzcd}
    \SMCATlcl \ar[rr, "S"] \ar[dr, hookrightarrow, "i"'] & & 
    \SMENR^{\lax,\cl}_{\CA} \ar[dl, "P"] \\
     & \SMCATl
  \end{tikzcd}
\]
where $i$ denotes the fully faithful inclusion of
$\SMCATlcl$ in $\SMCATl$. The functor $P$ is a $(0,1)$-fibration by
\cref{propn:SMENRlcl cocart}, and we also have a lax transformation
$\alpha \colon SU \to \id$ from above. However, we cannot directly
apply \cref{propn:extend adjn} to this data, as the functor $U$ does
\emph{not} generally take $P$-cocartesian morphisms to equivalences:
the cocartesian morphism from $\eC$ over $\phi \colon \uV \to \uW$ is
taken by $U$ to a functor given on mapping spaces by
\[ \uV(\bbone, \eC(x,y)) \to \uW(\bbone, \phi \eC(x,y)), \]
which is not necessarily an equivalence. We are therefore forced to
consider an enlargement of $\SMENRlc$:
\begin{defn}
  We take
  \[ \tQ := \SMENRlc \times_{\SMCATlcl} \AR(\SMCATlcl) \]
  to be the pullback of $U$ and $\ev_{1}$, and write $j \colon
  \SMENRlc \to \tQ$ for the inclusion and $\pi \colon \tQ \to
  \SMENRlc$ for the projection. Then we have a commutative triangle
  \[
    \begin{tikzcd}
     \SMENRlc \ar[rr, "j"] \ar[dr, "U"'] & & \tQ \ar[dl, "U'"] \\
      & \SMCATlcl,
    \end{tikzcd}
  \]
  where $U'$ is the functor given by $\ev_{0}$.  Note that we then
  have a canonical (strong) natural transformation
  \[\beta \colon \id \to j \pi.\] We further define
  \[P' := P \pi \colon \tQ \to \SMCATl\] and
  \[S' := jS \colon \SMCATlcl \to \tQ;\] then we have an equivalence
  \[ P'S' = P \pi j S \simeq P S \simeq i,\]
  so that we have a commutative triangle
  \[
    \begin{tikzcd}
     \SMCATlcl \ar[rr, "S'"] \ar[dr, "i"'] & & \tQ \ar[dl, "P'"] \\
      & \SMCATl.
    \end{tikzcd}
  \]
  Note that we also have $U'S' = U'jS = US \simeq \id$.
\end{defn}

\begin{construction}
  By \cite[Theorem 4.2.1]{decorated}, the functor \[U' \colon \tQ \to
  \SMCATlcl\] is obtained from $U$ by freely adding cartesian
  $1$-morphisms; these are the commuting squares whose image under $\pi$
  is an equivalence. Moreover, by \cite[Corollary 3.5.1]{decorated}, 
  the functor \[U'_{*} \colon \FUN(\tQ, \tQ)^{\lax} \to
    \FUN(\tQ, \SMCATlcl)^{\lax}\]
  has cartesian lifts of all \emph{strong} transformations, given by
  strong transformations whose components are $U'$-cartesian. Consider
  the composite 
  \[ S'U' \xto{S'U'\beta} S'U'j\pi = jSU\pi \xto{j \alpha \pi} j\pi; \]
  using the equivalence $U'S' \simeq \id$, this lies over
  \[ U' \xto{U'\beta} U\pi\]
  as $U \alpha \simeq \id_{U}$. Since $\beta$ is a $U'_{*}$-cartesian
  morphism, it follows that we have a unique factorization
  \[ S'U' \xto{\alpha'} \id \xto{\beta} j\pi\]
  where $U'\alpha' \simeq \id_{U'}$. 
\end{construction}

\begin{propn}\label{propn:P' fib}
 $P'$ is a $(0,1)$-fibration. A morphism in $\tQ$ is
  $P'$-cocartesian \IFF{} its image under $\pi$ is $P$-cocartesian and
  its image under $U'$ is an equivalence. A 2-morphism is
  $P'$-cartesian \IFF{} its image under $P'$ is $P$-cartesian
  (equivalently, its image under $U\pi$ is an equivalence) and its
  image under $U'$ is an equivalence.
\end{propn}
\begin{proof}
  By definition, $P'$ is the composite
  \[ \tQ \xto{\pi} \SMENRlcl \xto{P} \SMCATl, \] where we know $P$ is
  a $(0,1)$-fibration from \cref{propn:SMENRlcl cocart}. Moreover,
  $\pi$ is the pullback along $U$ of
  $\ev_{1} \colon \AR(\SMCATlcl) \to \SMCATlcl$, which has cocartesian
  lifts of all $1$-morphisms by \cite[Theorem 4.1.2]{decorated}; the
  $\ev_{1}$-cocartesian morphisms are precisely those that go to
  equivalences under $\ev_{0}$. Thus $\tQ$ has $\pi$-cocartesian lifts
  of all morphisms, and these are precisely those that go to
  equivalences under $U'$. In particular, $\tQ$ has $\pi$-cocartesian
  lifts of $P$-cocartesian morphisms; this means that it has
  $P'$-cocartesian morphisms, and these are precisely those that map
  to $P'$-cocartesian morphisms under $\pi$ and equivalences under
  $U'$.  The same argument also shows that $\tQ$ has $P'$-cartesian
  2-morphisms, because all $P$-cartesian 2-morphisms map to
  equivalences under $U$, so we only need $\pi$ to have cartesian
  lifts of \emph{invertible} 2-morphisms.
\end{proof}

\begin{thm}\label{thm:S' U' extend adjn}
  There is a relative adjunction
  \[
    \begin{tikzcd}
     \AR^{\oplax}_{\mathrm{s}=\cl}(\SMCATl) \ar[rr, yshift=4pt,
     "\bar{S}'"{above}, "\perp"{below}]
     \ar[dr, "\ev_{1}"'] & & \tQ \ar[dl, "P'"] \ar[ll,
     yshift=-5pt, "\bar{U}'"{below}] \\
      & \SMCATl
    \end{tikzcd}
  \]
  where
  \begin{itemize}
  \item $\AR^{\oplax}_{\mathrm{s}=\cl}(\SMCATl) := \SMCATlcl
    \times_{\SMCATl} \AR^{\oplax}(\SMCATl)$ is the full sub-\itcat{}
    of $\AR^{\oplax}(\SMCATl)$ on morphisms whose source is closed,
  \item   $\bar{S}'$ is the functor obtained by extending $S'$ over the
  free fibration on $i$,
\item  $\bar{U}'$ corresponds to the functor $U'$
  together with the lax transformation $P'\alpha'$, 
\item the unit of the adjunction is an equivalence $\id \simeq
  \bar{U'}\bar{S'}$ (so $\bar{S'}$ is fully faithful),
\item the counit of the adjunction $\bar{\alpha}' \colon \bar{S}'\bar{U}' \to \id_{\tQ}$ is
    obtained by factoring $\alpha'$ through a cocartesian morphism
    with respect to the $(0,1)$-fibration
    \[ P'_{*} \colon \FUN(\tQ,\tQ)^{\lax} \to \FUN(\tQ, \SMCATl)^{\lax}.\]
  \end{itemize}
\end{thm}
\begin{proof}
  We apply \cref{propn:extend adjn}. For this we need to verify that:
  \begin{enumerate}
  \item $U'$ takes $P'$-cocartesian morphisms and $P'$-cartesian
    2-morphisms to equivalences.
  \item The composite 
    $U' \isoto U'S'U' \xto{U' \alpha'} U'$ is equivalent to $\id_{U'}$.
  \item The lax naturality squares of $\alpha'$ all contain
    $P'$-cartesian 2-morphisms,
  \item The lax naturality squares of $\alpha'$ commute for
    $P'$-cocartesian morphisms,
  \item The composite $S' \isoto S'U'S' \xto{\alpha' S'} S'$
    is equivalent to $\id_{S'}$.
  \end{enumerate}

  Point (1) is part of \cref{propn:P' fib}, while we have by
  definition of $\alpha'$ that 
  $U'\alpha' \simeq \id_{U'}$, which proves (2). Since we also have
  $\pi \alpha' = \alpha$,  (3) then follows from
  \cref{propn:alpha props}(iii) and \cref{propn:P' fib}. To prove (4),
  it then suffices to show that $P'$ takes the naturality 2-morphism
  for a $P'$-cocartesian morphism to an equivalence; this follows from
  \cref{propn:alpha props}(iv) as
  $P'\alpha' = P \pi \alpha' = P \alpha$. Finally, (5) follows from 
  \cref{propn:alpha props}(ii) as we have $\alpha' j \simeq j \alpha$, so that
  \[ \alpha' S' \simeq \alpha j S \simeq j \alpha S. \qedhere \]
\end{proof}

\begin{remark}
  Note that if instead of $\tQ$ we considered the full free
  $(1,0)$-fibration
  \[ \SMENRlc \times_{\SMCATlcl} \ARopl(\SMCATlcl) \xto{\ev_{0}}
    \SMCATlcl\] on $U$, the resulting version of $\alpha'$ would
  \emph{not} have cartesian 2-morphisms in its naturality
  squares. This is the reason we need to use $\AR(\SMCATlcl)$ to be
  able to apply \cref{propn:extend adjn}.
\end{remark}
 
Since the unit of the adjunction above is an equivalence, the functor
$\bar{S}'$ is fully faithful. Our next goal is to identify its image,
for which we introduce the following terminology:
\begin{defn}
  A \emph{pinning} of a
  closed symmetric monoidal $\oV$-\icat{} $\eC$ is a closed symmetric
  monoidal \icat{} $\oB$ together with a strong symmetric
  monoidal functor $f \colon \oB \to U(\eC)$ that is essentially
  surjective and preserves internal Hom's. We write
  \[ \PSMENRlcl \subseteq \tQ\]
  for the full sub-\itcat{} of the fibre product over $U$ and
  $\ev_{1}$ whose objects are the \emph{pinned} closed symmetric
  monoidal $\oV$-\icats{}, \ie{} closed symmetric monoidal
  $\oV$-\icats{} equipped with a pinning.
\end{defn}

\begin{cor}
  The adjunction of \cref{thm:S' U' extend adjn} restricts to an
  equivalence
  \[ \bar{S}' : \ARoplcl \simeq \PSMENRlcl : \bar{U}'.\]
  This restricts further to an equivalence between $\SMENRlcl$ and the
  full sub-\itcat{} of $\ARoplcl$ on the morphisms $\phi \colon \uC
  \to \uV$ such that the canonical natural transformation $\uC(\bbone, \blank) \to
  \uV(\bbone, \phi(\blank))$ is an equivalence.
\end{cor}
\begin{proof}
  Since $\bar{S}'$ is fully faithful, for the first statement we need
  to verify that its image is precisely the full sub-\itcat{}
  $\PSMENRlcl$. It suffices to check that for $(\eC, f
  \colon \uC' \to U\eC)$, the induced map $g\colon \phi_{*}\bar{\uC}' \to
  \eC$ is an equivalence \IFF{} $f$ is a pinning, where $\phi =
  \eC(\bbone,f(\blank)) \colon U\eC \to \uV$.

  We have a factorization of $f$ as
  \[ \uC' \xto{r} U\phi_{*}\bar{\uC}' \xto{Ug} U\eC \]
  where \cref{propn:SMENRlcl cocart} implies that $r$ is essentially surjective, strong symmetric
  monoidal, and preserves internal Homs. It follows that these three
  properties must hold for $f$ whenever $g$ is an equivalence, so that
  the image of $\bar{S}'$ is contained in $\PSMENRlcl$.

  Conversely, $g$ is an
  equivalence
  \IFF{} $Ug$ is strong symmetric 
  monoidal and essentially surjective, and $g$ gives an equivalence
  \[ \eC(\bbone, f[x,y]) \to \eC(\bbone, [fx,fy])\simeq \eC(fx,fy) \]
  on all
  mapping objects in $\oV$. From the factorization above we see that
  $Ug$ has the first two properties when $f$ does, and the last
  condition certainly holds if $f$ preserves
  internal Homs. This shows that the image of $\bar{S}'$ is precisely
  $\PSMENRlcl$, as required.

  For the second restriction, we observe that $\phi \colon \oC \to
  \oV$ will correspond to an object of $\SMENRlcl$ precisely when the
  induced map $\oC \to U\phi_{*}\bar{\oC}$ is an equivalence. Since
  this is essentially surjective and strong symmetric monoidal, this
  happens \IFF{} it is an equivalence on mapping spaces, that is if
  \[ \oC(x,y) \to \oV(\bbone,\phi[x,y]) \]
  is an equivalence for all $x,y$. Since the source is equivalent to
  $\oC(\bbone,[x,y])$, this holds precisely when the canonical maps
  \[ \oC(\bbone,c) \to \oV(\bbone, \phi c)\]
  are equivalences for all $c \in \oC$.  
\end{proof}

Restricting to a single fibre, we get:
\begin{cornoproof}[\cref{thm:main V}]\label{cor:smenrcladjfixed}
  For $\uV$ symmetric monoidal, there is an equivalence
  \[ \SMCATlcl_{\upslash \oV} \simeq \PSMENRlcl(\oV) \]
  between the \icat{} of pinned closed symmetric monoidal
  $\uV$-\icats{} and that of lax symmetric monoidal functors $\phi
  \colon \oC \to \oV$ such that $\oC$ is closed. This restricts to an
  equivalence between $\SMENRlcl(\oV)$ and the full subcategory on
  functors $\phi$ such that
  \[ \oC(\bbone,c) \isoto \oV(\bbone, \phi c)\]
  for $c \in \oC$. 
\end{cornoproof}

We also note how the equivalence restricts to strong symmetric
monoidal functors, using the following notation:
\begin{notation}
  We write $\PSMENRcl(\oV)$ for the wide sub-\itcat{} of
  $\PSMENRlcl(\oV)$ containing those morphisms
  \[ (F,f) \colon (\eC, \oC' \xto{p} U\eC) \to (\eD, \oD' \xto{q} U\eD)\]
  such that $f \colon \oC' \to \oD'$ is strong symmetric monoidal.
\end{notation}

\begin{observation}
  For a morphism $(F,f)$ in
  $\PSMENRlcl(\oV)$ as above, if $f$ is strong symmetric monoidal, then so is
  $F$: Given objects $x,y \in \eC$ we can choose $x',y' \in \oC'$ such
  that $x \simeq px', y \simeq py'$ as $p$ is essentially surjective.
  To see that the structure map $Fx \otimes Fy \to F(x \otimes y)$ is
  invertible, we can consider the structure map
  \[ Fpx'  \otimes Fp y' \to F(px' \otimes py') \isoto Fp(x' \otimes y')\]
  for $Fp$, since $p$ is strong. This is an equivalence as $Fp \simeq
  qf$ and by assumption $q$ and $f$ are both strong. In particular,
  the \itcat{}
  $\SMENRcl(\oV)$ of symmetric monoidal $\oV$-\icats{} and strong
  symmetric monoidal functors is the intersection of  $\PSMENRcl(\oV)$
  and $\SMENRlcl(\oV)$. Note, however, that the converse is not true:
  For $(F,f)$ above we may have that $F$ is strong, but $f$ is not.
\end{observation}

\begin{cornoproof}
  The equivalence of \cref{cor:smenrcladjfixed} restricts to an
  equivalence between $\PSMENRcl(\oV)$ and the wide
  subcategory
  \[ \SMCATcl \times_{\SMCATlcl} \SMCATlcl_{\upslash \oV} \subseteq
    \SMCATlcl_{\upslash \oV}\]
  whose morphisms are the oplax triangles
      \[
    \begin{tikzcd}[column sep=small]
     \oC \ar[rr, "F"] \ar[dr, ""{name=A}] & & \oD \ar[dl] \\
     & \oV
     \arrow[from=A,to=1-3,Rightarrow,shorten <=4pt, shorten >=10pt]
    \end{tikzcd}
  \]
  where $F$ is strong symmetric monoidal (but the functors to $\oV$
  are still lax). Moreover, the smaller subcategory  $\SMENRcl(\oV)$
  of strong symmetric monoidal functors
  corresponds to that where we impose the further condition that
    \[ \oC(\bbone,c) \isoto \oV(\bbone, \phi c)\]
  for $c \in \oC$.
\end{cornoproof}

\appendix

\section{Fibrations of $(\infty,2)$-categories}\label{app:fibs}
In this appendix we first briefly review the definition of fibrations
of \itcats{} and recall the construction of free fibrations. We then
use these to give a construction of adjunctions we apply to
closed symmetric monoidal enriched \icats{} in \cref{sec:extend adjn}.
\begin{defn}
  For $p \colon \tE \to \tB$ a functor of \itcats{}, a morphism $f
  \colon x \to y$ in $\tE$ is \emph{$p$-cocartesian} if for all $z
  \in \tE$ the commutative square of \icats{}
  \[
    \begin{tikzcd}
     \tE(y,z) \ar[r, "f^{*}"] \ar[d] & \tE(x,z)  \ar[d] \\
     \tB(py,pz) \ar[r, "p(f)^{*}"'] & \tB(px,pz)
    \end{tikzcd}
  \]
  is a pullback. Dually, $f$ is \emph{$p$-cartesian} if $p^{\op}$ is
  $p$-cocartesian.
  A 2-morphism $\alpha \colon f \to g$ in $\tE$ is
  \emph{weakly $p$-(co)cartesian} if $\alpha$ is a (co)cartesian
  morphism for the functor of \icats{}
  \[ \tE(x,y) \to \tB(px,py),\]
  and \emph{$p$-(co)cartesian} if the whiskering $h \alpha k$ is
  weakly $p$-cocartesian for all 1-morphisms $k \colon x' \to x$, $h
  \colon y \to y'$.
\end{defn}

\begin{defn}\label{def:01fib}
  We say that $p \colon \tE \to \tB$ is a \emph{$(0,1)$-fibration} if
  $\tE$ has all $p$-cocartesian lifts of $1$-morphisms in $\tB$ and
  all $p$-cartesian lifts of $2$-morphisms. Dually, $p$ is a
  \emph{$(1,0)$-fibration} if it has all $p$-cartesian lifts of
  $1$-morphisms and $p$-cocartesian lifts of 2-morphisms.
\end{defn}

\begin{notation}
  Recall that we can equip $\CatIT$ with the \emph{Gray tensor
    product}, which we just denote $\otimes$.
  This preserves colimits in each variable, and
  so has adjoints such that
  \[ \Map(\tA \otimes \tB, \tC) \simeq \Map(\tB, \FUN(\tA,
    \tC)^{\oplax}) \simeq \Map(\tA, \FUN(\tB, \tC)^{\lax}).\]
  The morphisms in $\FUN(\tA,\tB)^{\plax}$ are \emph{(op)lax
    transformations}. We abbreviate $\AR^{\plax}(\tC) := \FUN([1],
  \tC)^{\plax}$.
\end{notation}

The following is a special case of \cite[Theorem 4.2.1]{decorated}:
\begin{thm}
  For a functor of \itcats{} $p \colon \tA \to \tB$, let $\bar{p}
  \colon \FFr(\tA)
  \to \tB$ denote the pullback
  \[ \FFr(\tA) := \tA \times_{\tB} \ARopl(\tB)\]
  via $p$ and $\ev_{0}$, with the map $\bar{p}$ induced by evaluation
  at $1$. Then
  \begin{enumerate}[(1)]
  \item $\bar{p}$ is a $(0,1)$-fibration. A morphism is
    $\bar{p}$-cocartesian \IFF{} its projection to $\tA$ is an
    equivalence and its projection to $\ARopl(\tB)$ is a commuting
    square; a 2-morphism is cartesian \IFF{} its projection to $\tA$
    is an equivalence.
  \item The inclusion $i \colon \tA \to \FFr(\tA)$, induced by
  restriction along $s_{0} \colon [1] \to [0]$, exhibits $\bar{p}$ as
  the free $(0,1)$-fibration on $p$: for $\tE \to \tB$ a
  $(0,1)$-fibration, the functor
  \[ \FUN^{(0,1)}_{/\tB}(\FFr(\tA), \tE) \to \FUN_{/\tB}(\tA, \tE)\]
  is an equivalence, where the source denotes the \itcat{} of functors
  that preserve cocartesian morphisms and cartesian 2-morphisms.
  \end{enumerate}
\end{thm}

\begin{observation}\label{obs:free fib canonical lax tr}
  Let $\phi \colon [1] \otimes [1] \to [1] \times [1] \to [1]$ be the
  functor that sends $(1,1)$ to $1$ and the remaining objects to
  $0$ --- that is, $\phi$ is the lax (but actually commuting) square
  \[
    \begin{tikzcd}
     0 \ar[r, "="] \ar[d, "="'] & 0  \ar[d] \\
     0 \ar[r] & 1
    \end{tikzcd}
  \]
  in $[1]$. Composition with $\phi$ induces a functor
  \[ \ARopl(\tB) \to \FUN([1] \otimes [1], \tB) \simeq
    \ARopl(\ARopl(\tB)),\]
  which is adjoint to a lax transformation
  \[ \Phi\colon [1] \otimes \ARopl(\tB) \to \ARopl(\tB)\]
  from $s_{0}\ev_{0}$ to $\id$.
  Here for an arrow $a \to b$, $\Phi$ gives the commutative square
  \[
    \begin{tikzcd}
     a \ar[r, "="] \ar[d, "="'] & a  \ar[d] \\
     a \ar[r] & b, 
    \end{tikzcd}
  \]
  while $\ev_{0}\Phi$ is the identity. This shows that $\Phi$ is a
  cocartesian morphism for the $(0,1)$-fibration
  \[ \ev_{1,*} \colon \FUN(\ARopl(\tB),\ARopl(\tB))^{\lax} \to
    \FUN(\ARopl(\tB), \tB)^{\lax}\]
  by \cite[Corollary
  3.5.1]{decorated}: its  cocartesian morphisms are the lax
  transformations whose components are $\ev_{1}$-cocartesian morphisms
  (meaning they are commutative squares that go to equivalences under
  $\ev_{0}$) and whose lax naturality squares contain
  $\ev_{1}$-cartesian 2-morphisms (meaning they go to equivalences
  under $\ev_{0}$). Note that the image $\ev_{1}\Phi$ is the lax
  natural transformation $[1] \otimes \ARopl(\tB) \to \tB$ from
  $\ev_{0}$ to $\ev_{1}$ that is adjoint to the identity of
  $\ARopl(\tB)$.

  Moreover, for $p \colon \tA \to \tB$ we can
  pull
  back $\Phi$ to a lax transformation
  \[ \Phi_{p} \colon [1] \otimes \FFr(\tA) \to \FFr(\tA)\]
  from $\id$ to $i\pi$, where $\pi$ is the projection to
  $\tA$. This is again a cocartesian morphism for the
  $(0,1)$-fibration
  \[ \bar{p}_{*} \colon \FUN(\FFr(\tA), \FFr(\tA))^{\lax} \to
    \FUN(\FFr(\tA), \tB)^{\lax}. \]
\end{observation}

\begin{observation}\label{obs:map to free fib}
  Suppose given a functor $p \colon \tA \to \tB$, a commutative triangle
  \[
    \begin{tikzcd}
     \tE \ar[rr, "\bar{G}"] \ar[dr, "q"'] & & \FFr(\tA) \ar[dl, "\bar{p}"] \\
      & \tB
    \end{tikzcd}
  \]
  into the free $(0,1)$-fibration on $p$ is uniquely determined by a
  functor $G \colon \tE \to \tA$ and a lax transformation $\alpha
  \colon pG \to q$. Moreover, if $q$ is a $(0,1)$-fibration, then
  $\bar{G}$ is a morphism of $(0,1)$-fibrations \IFF{}
  \begin{itemize}
  \item $G$ takes $p$-cocartesian morphisms and $p$-cartesian
    2-morphisms to equivalences in $\tA$,
  \item the lax naturality square of $\alpha$ commutes for a
    $p$-cocartesian morphism.
  \end{itemize}
\end{observation}

Our main use of free \itcatl{} fibrations will be via the following
construction of adjunctions:
\begin{propn}\label{propn:extend adjn}
  Given
  \begin{itemize}
  \item a $(0,1)$-fibration $p \colon \tE \to \tB$,
  \item functors $F \colon \tA \to \tE$ and $G \colon \tE \to \tA$,
  \item an equivalence $\id \isoto GF$,
  \item a lax transformation $\alpha \colon FG \to \id_{\tE}$,
  \end{itemize}
  such that
  \begin{enumerate}[(i)]
  \item $G$ takes $p$-cocartesian morphisms and $p$-cartesian
    2-morphisms to equivalences,
  \item the lax naturality squares of $\alpha$ all contain
    $p$-cartesian 2-morphisms and commute for $p$-cocartesian
    morphisms,
  \item the composites $F \isoto FGF \xto{\alpha F} F$ and
    $G \isoto GFG \xto{G \alpha} G$ are equivalent to the respective
    identities.
  \end{enumerate}
  Then there is an adjunction
  \[ \bar{F} : \FFr(\tA) \rightleftarrows \tE : \bar{G}\]
  over $\tB$ where
  \begin{itemize}
  \item $\bar{q} \colon \FFr(\tA) \to \tB$ is the free $(0,1)$-fibration on $q := pF$
    and $\bar{F}$ is the unique morphism of $(0,1)$-fibrations
    extending $F$,
  \item $\bar{G}$ is the functor determined by $G$ and $p\alpha$ as in
    \cref{obs:map to free fib}, and this is a morphism of $(0,1)$-fibrations,
  \item the unit of the adjunction is an equivalence $\id \isoto
    \bar{G}\bar{F}$ (so $\bar{F}$ is fully faithful),
  \item the counit of the adjunction $\bar{\alpha} \colon \bar{F}\bar{G} \to \id_{\tE}$ is
    obtained by factoring $\alpha$ through a cocartesian morphism over
    $p\alpha$ with respect to the $(0,1)$-fibration
    \[ p_{*} \colon \FUN(\tE,\tE)^{\lax} \to \FUN(\tE, \tB)^{\lax}.\]
  \end{itemize}
\end{propn}
\begin{proof}
  The conditions on $\alpha$ imply that $\bar{G}$ is a morphism of
  $(0,1)$-fibrations by \cref{obs:map to free fib}.

  Let $i$ denote the inclusion $\tA \to \FFr(A)$ and $\pi$ the
  projection $\FFr(\tA) \to \tA$. Then by \cref{obs:free fib canonical
    lax tr} we have a canonical lax natural transformation
  $\beta \colon i \pi \to \id_{\FFr(A)}$ that is a cocartesian morphism for the
  $(0,1)$-fibration
  \[ \bar{q}_{*} \colon
    \FUN(\FFr(\tA),\FFr(\tA))^{\lax} \to
    \FUN(\FFr(\tA), \tB)^{\lax}.\] Composing with $\bar{F}$
  and $\bar{G}$ and using that $\bar{F}i = F$ and $\pi \bar{G} = G$,
  we get lax transformations $\beta \bar{G} \colon iG \to \bar{G}$ and
  $\bar{F} \beta \colon F \pi \to \bar{F}$; since $\bar{F}$ and
  $\bar{G}$ are morphisms of $(0,1)$-fibrations, these are cocartesian for the $(0,1)$-fibrations
  \[  \FUN(\FFr(\tA), \tE)^{\lax} \xto{p_{*}} \FUN(\FFr(\tA),
    \tB)^{\lax}, \quad \FUN(\tE, \FFr(\tA))^{\lax} \xto{\bar{q}_{*}} \FUN(\tE, \tB)^{\lax}, \]
  lying over $p\alpha$ and $\bar{q}\beta$, respectively. Moreover, the
  composites
  \[ \gamma \colon iGF\pi \xto{\beta\bar{G} F\pi} \bar{G}F\pi
    \xto{\bar{G}\bar{F}\beta} \bar{G}\bar{F},\qquad
     \delta \colon FG \simeq F\pi i G \xto{\bar{F}\beta i G} \bar{F} i G
    \xto{\bar{F}\bar{G}\beta} \bar{F}\bar{G} 
  \]
  are also cocartesian for the corresponding lax functor \itcats{}.

  Here assumption (iii) implies that  the composite
  \[ i\pi \simeq iGF\pi  \to \bar{G}F \pi\]
  lies over the identity, so the composite of $\gamma$ with this
  equivalence lies over
  $\bar{q}
  \beta$. By the uniqueness of cocartesian morphisms we can then
  define the unit equivalence so that we get a
  commutative square
  \[
    \begin{tikzcd}
     i\pi  \ar[r, "\sim"] \ar[d, "\beta"] & iGF\pi  \ar[d, "\gamma"] \\
     \id \ar[r, "\sim"'] & \bar{G}\bar{F},
    \end{tikzcd}
  \]
  identifying the two vertical cocartesian morphisms over
  $\bar{q}\beta$. Note that composing this square with $\bar{F}$ and $\bar{G}$ now gives commutative
  squares
  \[
    \begin{tikzcd}
     F\pi  \ar[r, "\sim"] \ar[d] & FGF\pi  \ar[d] \\
     \bar{F} \ar[r, "\sim"'] & \bar{F}\bar{G}\bar{F},
   \end{tikzcd}
   \qquad
       \begin{tikzcd}
     iG  \ar[r, "\sim"] \ar[d] & iGFG  \ar[d] \\
     \bar{G} \ar[r, "\sim"'] & \bar{G}\bar{F}\bar{G},
   \end{tikzcd}
 \]
 where the vertical maps in both are cocartesian. 

  Our lax transformation $\alpha \colon FG \to \id$ factors through the cocartesian
  morphism over $p \alpha$ with source $FG$, which is $\delta$. Thus
  this factorization indeed yields a lax transformation $\bar{\alpha}
  \colon \bar{F}\bar{G} \to \id$. We claim this is in fact a strong
  natural transformation, \ie{} its a priori lax naturality squares
  actually commute. Indeed, the lax naturality squares of $\alpha$ are
  the composites of those for $\delta$ and $\bar{\alpha}$; since those
  for $\alpha$ by assumption contain cartesian 2-morphisms, as do
  those for $\delta$, the 3-for-2 property of cartesian morphisms
  implies that for $f \colon X \to Y$ in $\tE$, the composite of
  $\alpha_{f} \colon f \alpha_{X} \to \alpha_{Y} \bar{F}\bar{G} f$
  with the cocartesian morphism $FGX \to \bar{F}\bar{G}X$ is a
  cartesian 2-morphism. By \cite[Lemma 2.4.2]{decorated}, we can
  detect cartesian 2-morphisms after precomposition with cocartesian
  morphisms, so the naturality squares of $\bar{\alpha}$ must contain
  cartesian 2-morphisms; since they also lie over the identity, this
  shows they commute.

  Now the definition of $\bar{\alpha}$ means we also have commutative
  diagrams
 \[
    \begin{tikzcd}
     FGF\pi \ar[r, "\alpha F \pi"] \ar[d] & F \pi  \ar[d] \\
     FG \bar{F} \ar[r, "\alpha \bar{F}"'] \ar[d] & \bar{F} \ar[d, equals] \\
     \bar{F}\bar{G}\bar{F} \ar[r, "\bar{\alpha}\bar{F}"'] & \bar{F},
    \end{tikzcd}
    \qquad
    \begin{tikzcd}
     iGFG \ar[r, "iG\alpha"] \ar[d] & iG  \ar[d] \\
     \bar{G}FG \ar[r, "\bar{G}\alpha"'] \ar[d] & \bar{G} \ar[d, equals] \\
     \bar{G}\bar{F}\bar{G} \ar[r, "\bar{G}\bar{\alpha}"'] & \bar{G}
    \end{tikzcd}
  \]
  where the vertical maps are cocartesian. Combining these with the
  squares above we get diagrams
  \[
    \begin{tikzcd}
      F\pi \ar[r, "\sim"] \ar[d] & FGF\pi  \ar[d] \ar[r, "\alpha F
      \pi"] & F \pi \ar[d] \\
     \bar{F} \ar[r, "\sim"'] & \bar{F}\bar{G}\bar{F} \ar[r, "\bar{\alpha}\bar{F}"] & \bar{F},
   \end{tikzcd}
   \qquad
       \begin{tikzcd}
      iG \ar[r, "\sim"] \ar[d] & iGFG  \ar[d] \ar[r, "i G \alpha"] & IG \pi \ar[d] \\
     \bar{G} \ar[r, "\sim"'] & \bar{G}\bar{F}\bar{G} \ar[r, "\bar{G}\bar{\alpha}"] & \bar{G},
    \end{tikzcd}
  \]
  where in both cases the composite in the top row is equivalent to
  the identity and the vertical maps are cocartesian. By the
  uniqueness of cocartesian morphisms the composites in the bottom
  rows must therefore also be identities, which gives the adjunction
  identities we wanted to establish.
\end{proof}


\begin{bibdiv}
\begin{biblist}
\bib{decorated}{article}{
  author={Abell\'an, Fernando},
  author={Haugseng, Rune},
  author={Martini, Louis},
  title={Free fibrations, lax colimits and Kan extensions for $(\infty ,2)$-categories},
  date={2026},
  eprint={arXiv:2602.07604},
}

\bib{enropd}{article}{
  author={Chu, Hongyi},
  author={Haugseng, Rune},
  title={Enriched $\infty $-operads},
  eprint={arXiv:1707.08049},
  journal={Adv. Math.},
  volume={361},
  date={2020},
  pages={106913, 85},
}

\bib{GRV}{article}{
  author={Gaitsgory, Dennis},
  author={Rozenblyum, Nick},
  author={Varshavsky, Yakov},
  title={Applications of (higher) categorical trace I: the definition of AGCat},
  date={2025},
  eprint={arXiv:2512.21726},
}

\bib{enr}{article}{
  author={Gepner, David},
  author={Haugseng, Rune},
  title={Enriched $\infty $-categories via non-symmetric $\infty $-operads},
  journal={Adv. Math.},
  volume={279},
  pages={575--716},
  eprint={arXiv:1312.3178},
  date={2015},
}

\bib{nmorita}{article}{
  author={Haugseng, Rune},
  title={The higher {M}orita category of $E_{n}$-algebras},
  date={2017},
  eprint={arXiv:1412.8459},
  journal={Geom. Topol.},
  volume={21},
  number={3},
  pages={1631--1730},
}

\bib{HHLN1}{article}{
  author={Haugseng, Rune},
  author={Hebestreit, Fabian},
  author={Linskens, Sil},
  author={Nuiten, Joost},
  title={Lax monoidal adjunctions, two-variable fibrations and the calculus of mates},
  journal={Proc. Lond. Math. Soc. (3)},
  volume={127},
  date={2023},
  number={4},
  pages={889--957},
  eprint={arXiv:2011.08808},
}

\bib{unfolding}{article}{
  date={2026},
  author={Haugseng, Rune},
  author={Nikolaus, Thomas},
  title={Unfolding of symmetric monoidal $(\infty ,n)$-categories},
}

\bib{HeineEnr}{article}{
  author={Heine, Hadrian},
  title={An equivalence between enriched $\infty $-categories and $\infty $-categories with weak action},
  journal={Adv. Math.},
  volume={417},
  date={2023},
  pages={Paper No. 108941, 140},
  eprint={arXiv:2009.02428},
}

\bib{HeineBienr}{article}{
  author={Heine, Hadrian},
  title={On bi-enriched $\infty $-categories},
  date={2024},
  eprint={arXiv:2406.09832},
}

\bib{HeineColim}{article}{
  author={Heine, Hadrian},
  title={The higher algebra of weighted colimits},
  date={2024},
  eprint={arXiv:2406.08925},
}

\bib{HeyerMann}{article}{
  author={Heyer, Claudius},
  author={Mann, Lucas},
  title={6-Functor Formalisms and Smooth Representations},
  date={2024},
  eprint={arXiv:2410.13038},
}

\bib{HinichBVTens}{article}{
  author={Hinich, Vladimir},
  title={On diagrams of algebras},
  date={2023},
  eprint={arXiv:2304.05796},
}

\bib{KontsevichSoibelman}{article}{
  author={Kontsevich, Maxim},
  author={Soibelman, Yan},
  title={Deformations of algebras over operads and the Deligne conjecture},
  conference={ title={Conf\'erence Mosh\'e{} Flato 1999, Vol. I (Dijon)}, },
  book={ series={Math. Phys. Stud.}, volume={21}, publisher={Kluwer Acad. Publ., Dordrecht}, },
  date={2000},
  pages={255--307},
}

\bib{LurieCob}{article}{
  author={Lurie, Jacob},
  title={On the classification of topological field theories},
  conference={ title={Current developments in mathematics, 2008}, },
  book={ publisher={Int. Press, Somerville, MA}, },
  date={2009},
  pages={129--280},
  eprint={http://math.ias.edu/~lurie/papers/cobordism.pdf},
}

\bib{HA}{book}{
  author={Lurie, Jacob},
  title={Higher Algebra},
  date={2017},
  note={Available at \url {http://math.ias.edu/~lurie/}.},
}

\bib{McClureSmith}{article}{
  author={McClure, James E.},
  author={Smith, Jeffrey H.},
  title={A solution of Deligne's Hochschild cohomology conjecture},
  conference={ title={Recent progress in homotopy theory}, address={Baltimore, MD}, date={2000}, },
  book={ series={Contemp. Math.}, volume={293}, publisher={Amer. Math. Soc., Providence, RI}, },
  date={2002},
  pages={153--193},
  doi={10.1090/conm/293/04948},
}

\bib{ReutterZetto}{article}{
  author={Reutter, David},
  author={Zetto, Markus},
  title={Enriched $\infty $-categories as marked module categories},
  date={2025},
  eprint={arXiv:2501.07697},
}

\bib{Scholze6}{article}{
  author={Scholze, Peter},
  title={Six-functor formalisms},
  date={2025},
  eprint={https://people.mpim-bonn.mpg.de/scholze/SixFunctors.pdf},
}

\bib{ZettoOpd}{article}{
  author={Zetto, Markus},
  title={Enriched $\infty $-operads as marked algebras},
  date={2026},
  eprint={arXiv:2607.06676},
}
\end{biblist}
\end{bibdiv}
\end{document}